\documentclass[
	openany, 
	11pt, 
	a4paper, 
]{book} 

\usepackage[T2A]{fontenc}
\usepackage[utf8]{inputenc}
\usepackage[english]{babel}

\usepackage{amssymb}
\usepackage{amsmath}
\usepackage{amscd}
\usepackage{mathtools}
\usepackage{fancyhdr}
\usepackage{bm}
\usepackage{tikz}
\usepackage{graphicx}
\usepackage{wrapfig}
\usepackage{microtype}
\usepackage{physics}
\usepackage{hyperref}
\usepackage{caption}
\usepackage{rsfso}
\usepackage[mathscr]{eucal}

\usepackage{geometry}
\usepackage{fancyhdr} 
\newcommand{\Title}[1]{\begin{center}\baselineskip=6.0mm{\Large\textbf{#1}}\end{center}\vspace*{0.5mm}}
\newcommand{\Author}[1]{\centerline{\large\textbf{#1}}\vspace*{0mm}}

\newcommand{\SuperScript}[1]{\textsuperscript{\normalfont #1}}

\newcommand{\Star}{\textasteriskcentered{}}
\newcommand{\TwoStars}{\Star{}\textasteriskcentered{}}
\newcommand{\ThreeStars}{\TwoStars{}\textasteriskcentered{}}

\renewcommand{\tilde}{\widetilde}

\begin{document}
\selectlanguage{english}
\sloppy

\renewcommand{\abstractname}{}

\Title{A Mixed Problem for~the~Wave Equation in a Curvilinear Half-Strip with Discontinuous Initial Data}%

\Author{Viktor~I.~Korzyuk\SuperScript{1,2\Star}, Jan~V.~Rudzko\SuperScript{1\TwoStars}, Vladislav~V.~Kolyachko\SuperScript{1\ThreeStars}}

\begin{center}
\SuperScript{1}\textit{Institute of~Mathematics of~the~National Academy of~Sciences of~Belarus, Minsk, Belarus} \\ 
\SuperScript{2}\textit{Belarusian State University, Minsk, Belarus} \\

\vspace*{0.5mm}

\textit{e-mail: \SuperScript{\Star}korzyuk@bsu.by, \SuperScript{\TwoStars}janycz@yahoo.com, \SuperScript{\ThreeStars}vladislav.kolyachko@gmail.com} %
\let\thefootnote\relax\footnote{The~authors' research was supported by~the~Moscow Center for~Fundamental and Applied Mathematics of~M.~V.~Lomonosov Moscow State University under agreement № 075-15-2025-345.}
\end{center}

\noindent\makebox[\textwidth][c]{
\begin{minipage}{0.92\textwidth}
\small{
\textbf{Abstract.} 
For a~one-dimensional wave equation, we consider a~mixed problem in a curvilinear half-strip. The~initial conditions have a~first-kind discontinuity at one point. The~mixed problem models the~problem of~a longitudinal impact on a finite elastic rod with a movable boundary. We construct the~solution using the~method of~characteristics in an explicit analytical form. For the~problem in~question, we prove the~uniqueness of~the~solution and establish the~conditions under which its classical solution exists.
}
\end{minipage}
}

\bigskip

\textbf{1. Statement of~the~problem.} In the~curvilinear domain $Q=\{(t,x)~:\ t\in (0,\infty )\wedge x\in (\gamma (t),l)\}$, where $l$ is a~positive real number, of~two independent variables $(t,x)\in \overline{Q}\subset {{\mathbb{R}}^{2}}$, for~the~wave equation 
$$
(\partial _{t}^{2}-{{a}^{2}}\partial _{x}^{2})u(t,x)=f(t,x),\ \ \ \ (t,x)\in Q, \eqno(1)
$$
we consider the~following mixed problem with the~initial conditions
$$u(0,x)=\varphi (x),\,\,{{\partial }_{t}}u(0,x)=\psi (x)+\left\{ \begin{matrix}
   0, & x\in [0,l),  \\
   v, & x=l,  \\
\end{matrix} \right.\,\,x\in [0,l], \eqno(2)
$$
and the~boundary conditions
$$
u(t,\gamma (t))={{\mu }_{1}}(t), \quad (\partial _{t}^{2}+b{{\partial }_{x}})u(t,l)={{\mu }_{2}}(t), \quad t\in [0,\infty ), \eqno(3)
$$ 
where $a$, $v$, and $b$ are real numbers, $a>0$ for~definiteness, $f$ is a function given on the~set $\overline{Q}$, $\varphi$ and $\psi$ are some real-valued functions defined on the~segment $[0,l]$, and $\mu_1$ and $\mu_2$ are some real-valued functions defined on the~half-line $[0,\infty )$. We also assume that 
$$
\gamma \in {{C}^{1}}([0,\infty )),\ {\gamma }'(t)\in (-a,a) \text{ for~all } t\in [0,\infty ),\ \lim\limits_{t\to +\infty}\gamma (t)\pm at=\pm \infty, \eqno(4)
$$
and the~curves $x=\gamma (t)$ and $x=l$ do not intersect. 

The~mixed problem (1)--(4) models the~following problem from the~theory of~longitudinal impact~[1]. Suppose that an elastic finite homogeneous rod of~constant cross-section, whose left moving boundary $x=\gamma (t)$ is fixed, is subjected at the~initial moment $t=0$ to an impact at the~end $x=l$ by~a~load that sticks to the~rod. We also assume that an external volumetric force acts on the~rod, that the~displacements of~the~rod and the~rate of~their change at the~initial moment $t=0$ are not equal to zero, and that there are no shock waves in the~rod. Then, neglecting both the~weight of~the~rod as a~force and its possible vertical displacements, the~displacements $u(t,x)$ of~the~rod satisfy the~mixed problem (1)--(4), where $a=\sqrt{E{{\rho }^{-1}}}$, $b=SE{{M}^{-1}}$, where $E>0$ is Young's modulus of~the~rod material, $\rho >0$ is the~density of~the~rod material, $S>0$ is the~cross-sectional area~of~the~rod, $M>0$ is the~mass of~the~impacting load, $-v$ is the~velocity of~the~impacting load, ${{\mu }_{2}}$ is the~external force acting on the~end of~the~rod divided by~the~mass of~the~impacting load. The~quantity $-{{\mu }_{1}}(t)$ has a~physical meaning of~the~external force acting on the~end of~the~rod, ${{\mu }_{2}}(t)$ has a~physical meaning of~the~external force acting on the~end of~the~rod, divided by~the~mass of~the~impacting load. The~function f is the~external volumetric force divided by~$\rho$.

In the~case 
$$
\gamma (t)=0, \quad {{\mu }_{1}}={{\mu }_{2}}\equiv 0, \quad \varphi =\psi \equiv 0, \quad f\equiv 0, \eqno(5)
$$
J. Boussinesq~[2] constructed a formal solution to the~problem (1)--(4) using the~method of~characteristics. This approach was developed in [3--8]. E. L. Nikolai [4] found a general expression for~the~solution of~the~problem (1)--(5) in the~form of~a piecewise given function using the~Boussinesq method. S.~I.~Gaiduk~[9] solved the~problem (1)--(5) by~the~method of~contour integration~[10]. He strictly proved the~existence and uniqueness of~a~unique generalized solution, but not its physical correctness. A similar mixed problem, but with a boundary condition $(\partial _{t}^{2}+b{{\partial }_{x}}+c)u(t,l)=0$ instead of~$(\partial _{t}^{2}+b{{\partial }_{x}})u(t,l)=0$, was studied in the~work~[11] by~the~method of~contour integration, where again a~unique generalized solution was constructed and its physical correctness was not justified. Yufeng~X. and Dechao~Z. [12] obtained a formal analytical solution in the~form of~a trigonometric series to a problem that is similar to a mixed problem (1)--(5), but with the~boundary condition $u(t,0)-\beta {{\partial }_{x}}u(t,0)=0$ instead of~$u(t,0)=0$. The~problem (1)--(5) has also been solved using numerical methods, such as symbolic computations~[13] and the~finite element method [14].

When the~data is smooth, $v=0$, and the~half-strip is straight, i.e., $\gamma \equiv 0$, the~problem (1)--(4) has been studied using Fourier series [3, 15, 16] and the~method of~characteristics [17]. Auxiliary issues related to the~basis property of~the~system of~functions appearing in the~Fourier method for~the~problem (1)--(4) with $\gamma \equiv 0$ were studied in [18, 19]. Similar problems in curvilinear domains have been considered in the~works [20--22]. Questions related to the~stabilization and controllability of~solutions to the~wave equations in curvilinear domains have been studied in [23--25].

\textbf{2. Curvilinear half-strip.} Let us note some properties of~the~domain $Q$ in which the~problem is considered.

\textbf{Assertion 1.} {\it Let $({{t}_{0}},{{x}_{0}})\in Q$. Then the~value ${{x}_{0}}+a{{t}_{0}}$ is nonnegative under the~conditions $(4)$. }

\textbf{Assertion 2.} {\it Let $\alpha \in [0,\infty )$. Then the~equation $\gamma (t)+at=\alpha$ has a unique solution under the~conditions $(4)$.}

\textbf{Assertion 3.} {\it Let $\alpha \in (-\infty ,0]$. Then the~equation $\gamma (t)-at=\alpha$ has a~unique solution under the~conditions $(4)$.}

\textbf{Assertion 4.} {\it Let $({{t}_{0}},{{x}_{0}})\in Q$. Then the~curve $(t,\gamma (t))$ intersects the~line $x+at={{x}_{0}}+a{{t}_{0}}$ at~a~single point under the~conditions $(4)$.}

\textbf{Assertion 5.} {\it Let $({{t}_{0}},{{x}_{0}})\in Q$ and ${{x}_{0}}-a{{t}_{0}}\le 0$. Then the~curve $(t,\gamma (t))$ intersects the~line $x-at={{x}_{0}}-a{{t}_{0}}$ at a single point under the~conditions $(4)$.}

The~\textbf{proofs} of~Assertions~1 -- 5 are given in the~article~[21].

Consider the~following functions:
$$
{{\gamma }_{+}}:[0,\infty )\ni t\mapsto \gamma (t)+at, \quad {{\gamma }_{-}}:[0,\infty )\ni t\mapsto \gamma (t)-at
$$

We also need the~inverse of~the~functions ${{\gamma }_{+}}$ and ${{\gamma }_{-}}$, which will be denoted by~the~symbols ${{\Phi }_{+}}$ and ${{\Phi }_{-}}$, respectively, i.e., ${{\Phi }_{+}}(\gamma (t)+at)=t$ and ${{\Phi }_{-}}(\gamma (t)-at)=t$. Such functions exist by~Assertions~2 and~3. From the~inverse function theorem, we get the~formulas:
\begin{align}
    & {{{\Phi }'}_{-}}(t)=\frac{1}{{\gamma }'({{\Phi }_{-}}(t))-a}, \quad {{{\Phi }''}_{-}}(t)=-\frac{{\gamma }''({{\Phi }_{-}}(t))}{{{\left( {\gamma }'({{\Phi }_{-}}(t))-a \right)}^{3}}}, \quad t\in [0,\infty ), \tag{6} \\
    & {{{\Phi }'}_{+}}(t)=\frac{1}{{\gamma }'({{\Phi }_{+}}(t))+a}, \quad {{{\Phi }''}_{+}}(t)=-\frac{{\gamma }''({{\Phi }_{+}}(t))}{{{\left( {\gamma }'({{\Phi }_{+}}(t))+a \right)}^{3}}}, \quad t\in [0,\infty ) \tag{7}
\end{align}
Note that the~representations~(6) and~(7) and condition (4) imply that ${{\Phi }_{+}}$ is an increasing function and ${{\Phi }_{-}}$ is a decreasing function.

\textbf{3. Auxiliary problem.} Consider the~following simple case 
\[v=0. \eqno(8)\]
The~solution u of~the~problem (1)--(4), (8) has the~form
\[u(t,x)=w(t,x)+g(x-at)+p(x+at),\eqno(9)\]
where $w$ is a particular solution of~Eq.~(1). We can take it from the~paper~[20], it satisfies the~homogeneous initial conditions
\[w(0,x)={{\partial }_{t}}w(0,x)=0, \quad x\in [0,l],\]
and belongs to the~class ${{C}^{2}}(\overline{Q})$ if, for~example, $f\in {{C}^{1}}(\overline{Q})$. Moreover, $\partial _{t}^{2}w(0,x)=f(0,x)$ holds for~all $x\in [0,l]$.

So, we want to find closed-form expressions for~the~functions $g$ and $p$. To do this, we partition the~domain $\overline{Q}$ according to the~following formulas (for~clarity see Fig. 1):
\begin{equation}
    \begin{aligned}
        &{{Q}^{\left( 0,0 \right)}}=Q\cap \left\{ (t,x):x-at\in [0,l]\wedge x+at\in [0,l] \right\}\!,\\
&{{Q}^{\left( 1,0 \right)}}=Q\cap \left\{ (t,x):x-at\in [{{\gamma }_{-}}({{r}_{1}}),0]\wedge x+at\in [0,l] \right\}\!,\\
&{{Q}^{\left( 0,1 \right)}}=Q\cap \left\{ (t,x):x-at\in [0,l]\wedge x+at\in [l,l+a{{l}_{1}}] \right\}\!,\\
&{{Q}^{\left( i,j \right)}}=Q\cap \left\{ (t,x):x-at\in [{{\gamma }_{-}}({{r}_{i}}),{{\gamma }_{-}}({{r}_{i-1}})]\wedge x+at\in [l+a{{l}_{j-1}},l+a{{l}_{j}}] \right\}\!,
    \end{aligned}
    \tag{10}
\end{equation}
where ${{r}_{0}}={{l}_{0}}=0$, ${{l}_{i}}={{r}_{i-1}}+{{a}^{-1}}(l-\gamma ({{r}_{i-1}}))$, ${{r}_{i}}={{\Phi }_{+}}(l+a{{l}_{i-1}})$.

\begin{figure}[htb]
    \begin{center}
        \includegraphics[scale=0.76]{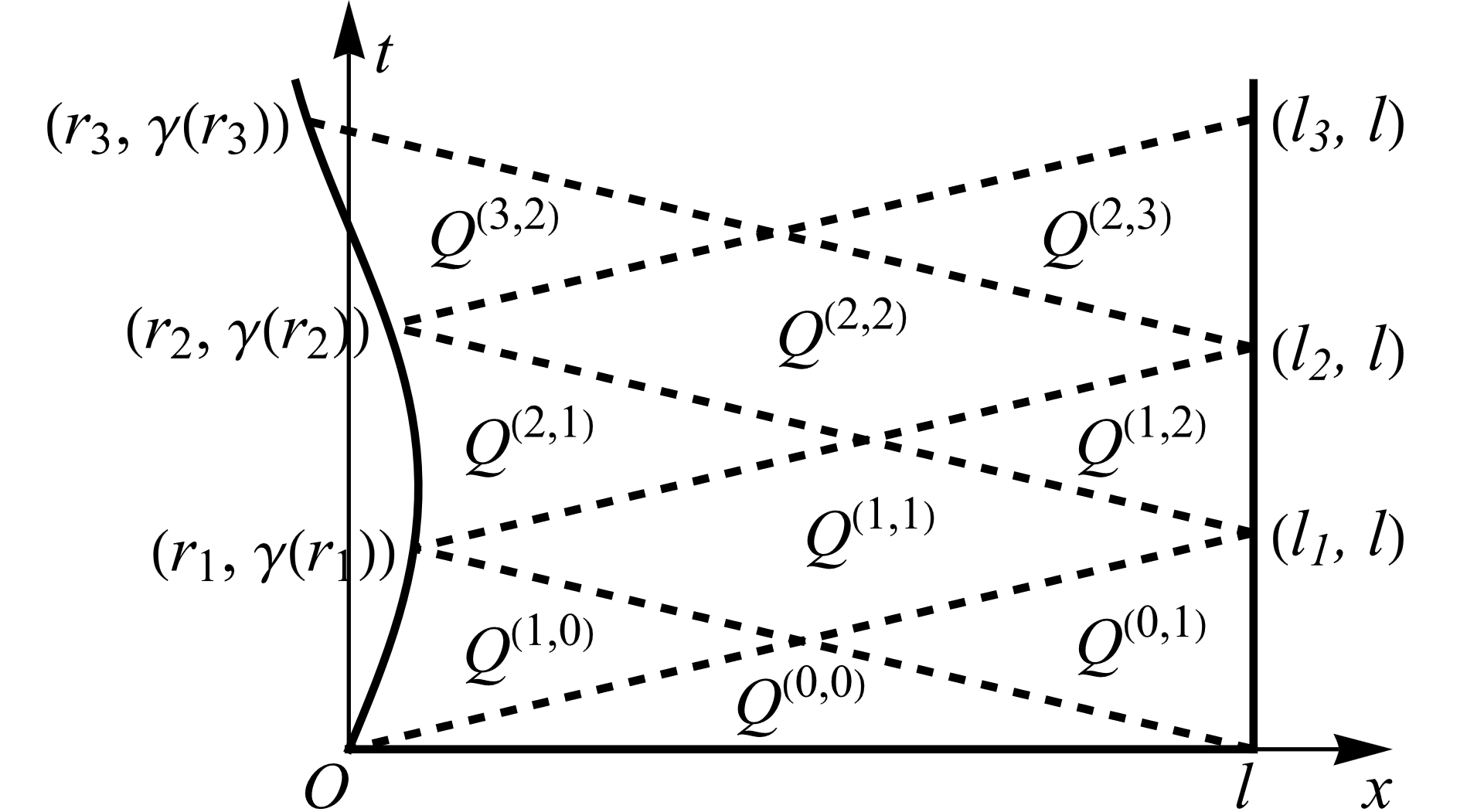}
    \end{center}
    \caption{Partitioning of~the~domain $Q$.}
\end{figure}

Let us demonstrate the~correctness of~the~partitioning (10) of~the~domain $Q$. First, we draw the~characteristic of~the~family $x-at=\text{const}$ through the~point $\left( {{r}_{0}}=0,\gamma\!\left( {{r}_{0}} \right)=0 \right)$, which intersects the~line $x=l$ at the~point $\left( {{l}_{1}}=l{{a}^{-1}},l \right)$. Next, we draw the~characteristic of~the~family $x+at=\text{const}$ through the~point $\left( {{l}_{1}}=l{{a}^{-1}},l \right)$. According to Assertion 4, this characteristic intersects the~curve $x=\gamma\!\left( t \right)$ at the~point $\left( {{r}_{2}},\gamma\!\left( {{r}_{2}} \right) \right).$ At the~same time, ${{r}_{2}}>0$. Next, we draw the~characteristic of~the~family $x-at=\text{const}$ through the~point $\left( {{r}_{2}},\gamma\!\left( {{r}_{2}} \right) \right)$. This characteristic intersects the~line $x=l$ at the~point $\left( {{l}_{3}},l \right)$ and, by~Assertion 5, does not intersect the~curve $x=\gamma\!\left( t \right)$ at any point other than $\left( {{r}_{2}},\gamma\!\left( {{r}_{2}} \right) \right)$. Thus, starting from $i=0$, we sequentially define the~points: $\left( {{l}_{i}},l \right)$, $\left( {{r}_{i+1}},\gamma\!\left( {{r}_{i+1}} \right) \right)$, $\left( {{l}_{i+2}},l \right)$, $\left( {{r}_{i+3}},\gamma\!\left( {{r}_{i+3}} \right) \right)$, $\left( {{l}_{i+4}},l \right)$, $\left( {{r}_{i+5}},\gamma\!\left( {{r}_{i+5}} \right) \right)$ etc. 

Let us describe the~process of~finding points of~the~form 
$$
\left( {{l}_{2i}},l \right),\quad \left( {{r}_{2i+1}},\gamma\!\left( {{r}_{2i+1}} \right) \right)\!, \quad i=0,1,\ldots 
$$
Through the~point $\left( {{l}_{0}}=0,l \right)$ we draw the~characteristic of~the~family $x+at=\text{const}$, which intersects the~curve $x=\gamma\!\left( t \right)$ at one point $\left( {{r}_{1}},\gamma\!\left( {{r}_{1}} \right) \right)$ by~virtue of~Assertion 5. Now, through the~point  $\left( {{r}_{1}},\gamma\!\left( {{r}_{1}} \right) \right)$ we draw the~characteristic of~the~family $x-at=\text{const}$, which intersects the~line $x=l$ at one point $\left( {{l}_{2}},l \right)$ and, by~Assertion 5, it does not intersect the~curve $x=\gamma\!\left( t \right)$ at any point other than $\left( {{r}_{1}},\gamma\!\left( {{r}_{1}} \right) \right)$ according to Assertion 5. Then through the~point $\left( {{l}_{2}},l \right)$ we draw the~characteristic of~the~family $x+at=\text{const}$, which intersects the~curve $x=\gamma\!\left( t \right)$ at one point $\left( {{r}_{3}},\gamma\!\left( {{r}_{3}} \right) \right)$ by~virtue of~Assertion 4. Thus, starting from $i=0$, we sequentially define the~points: $\left( {{r}_{i}},\gamma\!\left( {{r}_{i}} \right) \right)$, $\left( {{l}_{i+1}},l \right)$, $\left( {{r}_{i+2}},\gamma\!\left( {{r}_{i+3}} \right) \right)$, $\left( {{l}_{i+3}},l \right)$, $\left( {{r}_{i+4}},\gamma\!\left( {{r}_{i+4}} \right) \right)$ etc.

We will look for~the~function as a piecewise-defined, i.e.,
\[u(t,x)={{u}^{(i,j)}}(t,x)\quad (t,x)\in \overline{{{Q}^{(i,j)}}}. \eqno(11)\]
Because of~(9)--(11), we can write
\[{{u}^{(i,j)}}(t,x)=w(t,x)+{{g}^{(i)}}(x-at)+{{p}^{(j)}}(x+at),\quad (t,x)\in \overline{{{Q}^{(i,j)}}}. \eqno(12)\]

We determine the~functions ${{g}^{\left( 0 \right)}}$ and ${{p}^{\left( 0 \right)}}$ from the~Cauchy conditions~(2):
\begin{align}
{{g}^{\left( 0 \right)}}(x)=\frac{\varphi (x)}{2}-\frac{1}{2a}\int\limits_{0}^{x}{\psi (\xi )\, d\xi }+{{C}_{1}}, \quad x\in [0,l], \tag{13} \\
{{p}^{\left( 0 \right)}}(x)=\frac{\varphi (x)}{2}+\frac{1}{2a}\int\limits_{0}^{x}{\psi (\xi )\, d\xi }-{{C}_{1}}, \quad x\in [0,l], \tag{14}
\end{align}
where ${{C}_{1}}$ is a real number. The~function ${{g}^{\left( i \right)}}$ for~$i\in \mathbb{N}$ can be defined from the~Dirichlet boundary condition (3) on the~curve $x=\gamma (t)$. We substitute (12), where ${{Q}^{(i,j)}}={{Q}^{(i,i-1)}}$, into~(3) and obtain
\[
w(t,\gamma (t))+{{g}^{(i)}}(\gamma (t)-at)+{{p}^{(i-1)}}(\gamma (t)+at)={{\mu }_{1}}(t), \quad t\in [{{r}_{i-1}},{{r}_{i}}], \quad i\in \mathbb{N}.
\]
Changing the~variable $z=\gamma (t)-at$, i.e., $t={{\Phi }_{-}}(z)$, results in the~equation
\[w({{\Phi }_{-}}(z),\gamma ({{\Phi }_{-}}(z)))+{{g}^{(i)}}(z)+{{p}^{(i-1)}}(\gamma ({{\Phi }_{-}}(z))+a{{\Phi }_{-}}(z))={{\mu }_{1}}({{\Phi }_{-}}(z)),  \quad {{\Phi }_{-}}(z)\in [{{r}_{i-1}},{{r}_{i}}], \quad i\in \mathbb{N},\]
which we can solve to obtain
\[{{g}^{(i)}}(z)={{\mu }_{1}}({{\Phi }_{-}}(z))-{{p}^{(i-1)}}(\gamma ({{\Phi }_{-}}(z))+a{{\Phi }_{-}}(z))-w({{\Phi }_{-}}(z),\gamma ({{\Phi }_{-}}(z))),  \quad {{\Phi }_{-}}(z)\in [{{r}_{i-1}},{{r}_{i}}], \quad i\in \mathbb{N}. \eqno(15)\]     
The~function ${{p}^{\left( j \right)}}$ for~$j\in \mathbb{N}$ can be defined from the~boundary condition (3) on the~line $x=l$. Again, we substitute (12), where ${{Q}^{(i,j)}}={{Q}^{(j-1,j)}}$, into (3) and get
\begin{multline}
    {{a}^{2}}{{D}^{2}}{{g}^{(j-1)}}(l-at)+{{a}^{2}}{{D}^{2}}{{p}^{(j)}}(l+at)+b(D{{g}^{(j-1)}}(l-at)+D{{p}^{(j)}}(l+at)+{{\partial }_{x}}w(t,l))+\partial _{t}^{2}w(t,l)={{\mu }_{2}}(t), \\
    t\in [{{l}_{i-1}},{{l}_{i}}],  \quad j\in \mathbb{N},
    \tag{16}
\end{multline}
where $D$ is the~Newton--Leibniz operator. Changing the~variable $t={{a}^{-1}}(z-l)$ transforms Eq. (16) into
\begin{multline}
    {{a}^{2}}{{D}^{2}}{{p}^{\left( j \right)}}(z)+bD{{p}^{\left( j \right)}}(z)={{\mu }_{2}}\!\left( \frac{z-l}{a} \right)-{{a}^{2}}{{D}^{2}}{{g}^{\left( j-1 \right)}}(2l-z)-bD{{g}^{\left( j-1 \right)}}(2l-z) \, - \\ -b{{\partial }_{x}}v\!\left( \frac{z-l}{a},l \right)-\partial _{t}^{2}v\!\left( \frac{z-l}{a},l \right)\!, \quad
    z\in [l+a{{l}_{i-1}},l+a{{l}_{i}}], \quad j\in \mathbb{N}. \tag{17}
\end{multline}
We solve Eq. (17) and obtain
\begin{multline}
    {{p}^{\left( j \right)}}\!\left( z \right)={{p}^{\left( j \right)}}\!\left( l+a{{l}_{j-1}}+0 \right)+ \int\limits_{l+a{{l}_{j-1}}}^z \exp\!\left( \frac{b\left( l+a{{l}_{j-1}}-\eta  \right)}{{{a}^{2}}} \right)\bigg( D{{p}^{\left( j \right)}}\left( l+a{{l}_{j-1}}+0 \right) + \\ +\int\limits_{l+a{{l}_{j-1}}}^\eta {{a}^{-2}}\exp\!\left( \frac{b\left( \xi -l-a{{l}_{i-1}} \right)}{{{a}^{2}}} \right)
      \left\{ {{\mu }_{2}}\!\left( \frac{\xi -l}{a} \right)-{{a}^{2}}{{D}^{2}}{{g}^{\left( j-1 \right)}}\left( 2l-\xi  \right)-bD{{g}^{\left( j-1 \right)}}\left( 2l-\xi  \right)- \right. \\
    \left. -\left. b{{\partial }_{x}}v\!\left( \frac{\xi -l}{a},l \right)-\partial _{t}^{2}v\!\left( \frac{\xi -l}{a},l \right) \right\}d\xi  \right)d\eta,  \quad z\in [l+a{{l}_{j-1}},l+a{{l}_{j}}], \quad j\in \mathbb{N}.  \tag{18}
\end{multline}
We choose the~values ${{p}^{\left( j \right)}}\!\left( l+a{{l}_{j-1}}+0 \right)$ and $D{{p}^{\left( j \right)}}\!\left( l+a{{l}_{j-1}}+0 \right)$ in the~representation~(18) by~continuity, i.e.,
$$
{{p}^{\left( j \right)}}\!\left( l+a{{l}_{j-1}}+0 \right)={{p}^{\left( j-1 \right)}}\!\left( l+a{{l}_{j-1}}-0 \right)\!, \:\: D{{p}^{\left( j \right)}}\left( l+a{{l}_{j-1}}+0 \right)=D{{p}^{\left( j-1 \right)}}\!\left( l+a{{l}_{j-1}}-0 \right)\!, \:\: j\in \mathbb{N}. \eqno(19) 
$$
According to formulas (15), (18), and (19), the~following relations hold for~all $i\in \mathbb{N}$:

\begin{align}
    & \Delta _{g}^{i}={{g}^{\left( i+1 \right)}}\!\left( {{\gamma }_{-}}\!\left( {{r}_{i}} \right) \right)-{{g}^{\left( i \right)}}\!\left( {{\gamma }_{-}}\!\left( {{r}_{i}} \right) \right)={{p}^{\left( i-1 \right)}}\!\left( {{\gamma }_{+}}\!\left( {{r}_{i}} \right) \right)-{{p}^{\left( i \right)}}\!\left( {{\gamma }_{+}}\!\left( {{r}_{i}} \right) \right)=0, \notag \\
    & \Delta _{g}^{i}={{g}^{\left( i+1 \right)}}\!\left( {{\gamma }_{-}}\!\left( {{r}_{i}} \right) \right)-{{g}^{\left( i \right)}}\!\left( {{\gamma }_{-}}\!\left( {{r}_{i}} \right) \right)={{p}^{\left( i-1 \right)}}\!\left( {{\gamma }_{+}}\left( {{r}_{i}} \right) \right)-{{p}^{\left( i \right)}}\!\left( {{\gamma }_{+}}\!\left( {{r}_{i}} \right) \right)=0, \notag \\
    & \widetilde{\Delta }_{g}^{i}=D{{g}^{\left( i+1 \right)}}\!\left( {{\gamma }_{-}}\left( {{r}_{i}} \right) \right)-D{{g}^{\left( i \right)}}\!\left( {{\gamma }_{-}}\!\left( {{r}_{i}} \right) \right)=-\dfrac{a+{\gamma }'\!\left( {{r}_{i}} \right)}{a-{\gamma }'\!\left( {{r}_{i}} \right)}\left( D{{p}^{\left( i-1 \right)}}\left( {{\gamma }_{+}}\!\left( {{r}_{i}} \right) \right)-D{{p}^{\left( i \right)}}\left( {{\gamma }_{+}}\!\left( {{r}_{i}} \right) \right) \right)=0, \notag \\
    & \widetilde{\widetilde{\Delta }}_{g}^{i}={{D}^{2}}{{g}^{\left( i+1 \right)}}\!\left( {{\gamma }_{-}}\left( {{r}_{i}} \right) \right)-{{D}^{2}}{{g}^{\left( i \right)}}\!\left( {{\gamma }_{-}}\left( {{r}_{i}} \right) \right)=\dfrac{{{\left( a+{\gamma }'\!\left( {{r}_{i}} \right) \right)}^{2}}}{{{\left( a-{\gamma }'\!\left( {{r}_{i}} \right) \right)}^{2}}}\left( {{D}^{2}}{{p}^{\left( i-1 \right)}}\left( {{\gamma }_{+}}\!\left( {{r}_{i}} \right) \right)-{{D}^{2}}{{p}^{\left( i \right)}}\!\left( {{\gamma }_{+}}\!\left( {{r}_{i}} \right) \right) \right)\!,\notag \\
    &\Delta _{p}^{i}={{p}^{\left( i+1 \right)}}\!\left( l+a{{l}_{i}} \right)-{{p}^{\left( i \right)}}\!\left( l+a{{l}_{i}} \right)=0, \quad \widetilde{\Delta }_{p}^{i}=D{{p}^{\left( i+1 \right)}}\!\left( l+a{{l}_{i}} \right)-D{{p}^{\left( i \right)}}\!\left( l+a{{l}_{i}} \right)=0, \notag\\
    & \widetilde{\widetilde{\Delta }}_{p}^{i}={{D}^{2}}{{p}^{\left( i+1 \right)}}\!\left( l+a{{l}_{i}} \right)-{{D}^{2}}{{p}^{\left( i \right)}}\!\left( l+a{{l}_{i}} \right)=-{{D}^{2}}{{g}^{\left( i \right)}}\!\left( l-a{{l}_{i}} \right)+{{D}^{2}}{{g}^{\left( i-1 \right)}}\!\left( l-a{{l}_{i}} \right)- \notag \\
    &\quad\quad-\frac{b}{{{a}^{2}}}D{{g}^{\left( i \right)}}\!\left( l-a{{l}_{i}} \right)+\frac{b}{{{a}^{2}}}D{{g}^{\left( i-1 \right)}}\!\left( l-a{{l}_{i}} \right)\!. \tag{20}
\end{align}
By virtue of~the~expressions ${{l}_{i}}={{r}_{i-1}}+{{a}^{-1}}(l-\gamma ({{r}_{i-1}}))$ and ${{r}_{i}}={{\Phi }_{+}}(l+a{{l}_{i-1}})$ we have
$$
\Delta _{g}^{i}=\widetilde{\Delta }_{g}^{i}=\Delta _{p}^{i}=\widetilde{\Delta }_{p}^{i-1}=0, \quad \widetilde{\widetilde{\Delta }}_{g}^{i}=-\frac{{{\left( a+{\gamma }'\!\left( {{r}_{i}} \right) \right)}^{2}}}{{{\left( a-{\gamma }'\!\left( {{r}_{i}} \right) \right)}^{2}}}\widetilde{\widetilde{\Delta }}_{p}^{i-1}, \quad \widetilde{\widetilde{\Delta }}_{p}^{i}=-\widetilde{\widetilde{\Delta }}_{g}^{i-1}-\frac{b}{{{a}^{2}}}\widetilde{\Delta }_{g}^{i-1},\quad i\in \mathbb{N}, \eqno(21) 
$$
The~base of~the~recurrence relations (21) can be computed using the~representations (13), (14), (16), (18), and (19). So, after some simple calculations, we get
\begin{align}
    \Delta _{g}^{0}&={{\delta }_{0}}={{\mu }_{1}}( 0 )-\varphi ( 0 ),\quad \widetilde{\Delta }_{g}^{0}={{\delta }_{1}}=\frac{\psi \left( 0 \right)+{\gamma }'( 0 ){\varphi }'( 0 )-\mu_1'( 0 )}{a-{\gamma }'( 0 )}, \notag \\
    \widetilde{{\widetilde{\Delta }}}_{g}^{i}&={{\delta }_{2}}=\frac{1}{{{\left( a-{\gamma }'\!\left( 0 \right) \right)}^{3}}}\left( \left( \mu_1'\!\left( 0 \right)-\psi\!\left( 0 \right) \right){\gamma }''\!\left( 0 \right)-{{a}^{3}}{\varphi }''\!\left( 0 \right)+{{a}^{2}}{\gamma }'\!\left( 0 \right){\varphi }''\!\left( 0 \right)+{\gamma }'\!\left( 0 \right) \right.\times \notag \\
&\times \left( f\!\left( 0,0 \right)+2{\gamma }'\!\left( 0 \right){\psi }'\left( 0 \right)+({\gamma }'{\left( 0 \right)})^{2} {\varphi }''\!\left( 0 \right)-\mu_1''\!\left( 0 \right) \right)- \notag \\
&-a\! \left( f\!\left( 0,0 \right)+2{\gamma }'\!\left( 0 \right){\psi }'\!\left( 0 \right)+{\varphi }'\!\left( 0 \right){\gamma }''\!\left( 0 \right)+({\gamma }'{{\left( 0 \right)}})^{2}{\varphi }''\!\left( 0 \right)-\mu_1''\!\left( 0 \right) \right)\!, \notag \\
\Delta _{p}^{0}&={{\rho }_{0}}=0,\quad \widetilde{\Delta }_{p}^{0}={{\rho }_{1}}=0,\quad \widetilde{\widetilde{\Delta }}_{p}^{0}={{\rho }_{2}}=-\frac{f\!\left( 0,l \right)-{{\mu }_{2}}\!\left( 0 \right)+b{\varphi }'\!\left( 0 \right)+{{a}^{2}}{\varphi }''\!\left( l \right)}{{{a}^{2}}}. \tag{22}
\end{align}

The~following assertion holds. 

\textbf{Assertion 6.} {\it Let the~smoothness conditions 
\begin{equation}
    \varphi \in {{C}^{2}}([0,l]), \:\: \psi \in {{C}^{1}}([0,l]), \:\: {{\mu }_{1}}\in {{C}^{2}}([0,\infty )), \:\: {{\mu }_{2}}\in C([0,\infty )), \:\: \gamma \in {{C}^{2}}([0,\infty )), \:\: f\in {{C}^{1}}(\overline{Q})
\tag{23}
\end{equation}
be satisfied. Then the~functions g and p, defined by~the~formulas (13)--(15), (18), (19), and
\begin{equation}
\begin{gathered}
    g(z)={{g}^{(0)}}(z), \quad z\in [0,l], \quad g(z)={{g}^{(i)}}(z), \quad {{\Phi }_{-}}(z)\in [{{r}_{i-1}},{{r}_{i}}], \quad i\in \mathbb{N},\\
    p(z)={{p}^{(0)}}(z), \quad z\in [0,l], \quad p(z)={{p}^{(i)}}(z), \quad z\in [l+a{{l}_{j-1}},l+a{{l}_{j}}], \quad j\in \mathbb{N}.
\end{gathered} \tag{24}
\end{equation}
are twice continuously differentiable if and only if the~following matching conditions are satisfied 
\begin{align}
    & {{\mu }_{1}}\!\left( 0 \right)-\varphi \!\left( 0 \right)=0, \tag{25} \\
    & \mu_1'\!\left( 0 \right)-\psi \!\left( 0 \right)+{\gamma }'\!\left( 0 \right){\varphi }'\!\left( 0 \right)=0  \tag{26}\\
    & \mu_1''\!\left( 0 \right)-\left( {{a}^{2}}+{{\left( {\gamma }'\!\left( 0 \right) \right)}^{2}} \right){\varphi }''\!\left( 0 \right)-f\!\left( 0,0 \right)-2{\gamma }'\!\left( 0 \right){\psi }'\left( 0 \right)-{\gamma }''\!\left( 0 \right){\varphi }'\!\left( 0 \right)=0,  \tag{27} \\
    & {{\mu }_{2}}\!\left( 0 \right)-f\!\left( 0,l \right)-b{\varphi }'\!\left( 0 \right)-{{a}^{2}}{\varphi }''\!\left( l \right)=0. \tag{28} 
\end{align}
}

The~\textbf{proof} is based on the~formulas (13)--(15), (18), (19), (21), and (22). Using the~method of~mathematical induction, we will first demonstrate that
$$
{{g}^{\left( i \right)}}\in {{C}^{2}}\!\left( \mathfrak{D}\!\left( {{g}^{\left( i \right)}} \right) \right)\!,~~~~{{p}^{\left( i \right)}}\in {{C}^{2}}\!\left( \mathfrak{D}\!\left( {{p}^{\left( i \right)}} \right) \right)\!,~~~~i=0,1,\ldots 
$$
Indeed, if $\varphi \in {{C}^{2}}\!\left( \left[ 0,l \right] \right)$, $\psi \in {{C}^{1}}\!\left( \left[ 0,l \right] \right)$, then, according to the~expressions (13) and (14), ${{g}^{\left( 0 \right)}}\in {{C}^{2}}\!\left( \left[ 0,l \right] \right)$ and ${{p}^{\left( 0 \right)}}\in {{C}^{2}}\!\left( \left[ 0,l \right] \right)$. Thus, we have proved the~base case of~induction. Now suppose that 
$$
{{g}^{\left( i \right)}}\in {{C}^{2}}\!\left( \mathfrak{D}\!\left( {{g}^{\left( i \right)}} \right) \right)\!, \quad {{p}^{\left( i \right)}}\in {{C}^{2}}\!\left( \mathfrak{D}\!\left( {{p}^{\left( i \right)}} \right) \right)
$$
is true for~some nonnegative integer $i$. In this case, if the~conditions (23) are true, then according to the~representations (15) and (18) we have 
$$
{{g}^{\left( i+1 \right)}}\in {{C}^{2}}\!\left( \left[ {{\gamma }_{-}}\left( {{r}_{i}} \right),{{\gamma }_{-}}\left( {{r}_{i+1}} \right) \right] \right)\!, \quad {{p}^{\left( i+1 \right)}}\in {{C}^{2}}\!\left( \left[ l+a{{l}_{i}},l+a{{l}_{i+1}} \right] \right)
$$
Therefore, we have proved the~induction step. The~function $g$ is twice continuously differentiable everywhere except perhaps at the~points ${{\gamma }_{-}}\!\left( {{r}_{i}} \right)$, $i=0,1,\ldots $, and the~function $p$ is twice continuously differentiable everywhere except perhaps at the~points $l+a{{l}_{i}}$, $i=0,1,\ldots $. For the~functions $g$ and $p$ to be twice continuously differentiable everywhere in the~domain of~definition, it is necessary and sufficient that the~following conditions be satisfied:
$$
\Delta _{g}^{i}=\tilde{\Delta }_{g}^{i}=\widetilde{{\tilde{\Delta }}}_{g}^{i}=\Delta _{p}^{i}=\tilde{\Delta }_{p}^{i}=\widetilde{{\tilde{\Delta }}}_{p}^{i}=0, \quad i=0,1,\ldots \eqno(29)
$$
According to the~discontinuity representations (21) and (22), the~conditions (29) are true if and only if the~matching conditions (25)--(28) are satisfied. The~assertion is proved.

\textbf{Assertion 7.} {\it The~functions $g$ and $p$ defined by~the~formulas (13)--(15), (18), (19), and (24) are of~the~form $g=\widetilde{g}+{{C}_{1}}$, $p=\widetilde{p}-{{C}_{1}}$, where $\tilde{g}$ and $\tilde{p}$ are some functions not depending on the~constant $C_1$.}

\textbf{Proof.} We will prove the~theorem using the~method described in~[26, p.~179]. We will apply the~method of~mathematical induction. According to the~formulas (13)--(15), the~expression holds 
$$
{{g}^{\left( 0 \right)}}( x )={{\widetilde{g}}^{\left( 0 \right)}}( x )+{{C}_{1}}, \quad {{p}^{\left( 0 \right)}}( x )={{\widetilde{p}}^{\left( 0 \right)}}( x )-{{C}_{1}}, \quad x\in [ 0,l ],
$$
where
\begin{align*}
    &{{\widetilde{g}}^{\left( 0 \right)}}( x )=\frac{\varphi \left( x \right)}{2}-\frac{1}{2a}\int\limits_{0}^{x}{\psi \!\left( \xi  \right) d\xi },\ \ \ \ x\in [ 0,l ], \\
    &{{\widetilde{p}}^{\left( 0 \right)}}( x )=\frac{\varphi \left( x \right)}{2}+\frac{1}{2a}\int\limits_{0}^{x}{\psi \!\left( \xi  \right) d\xi },\ \ \ \ x\in [ 0,l ].
\end{align*}
Thus, the~base case of~induction is proved. Now, let us assume that the~relation 
$$
{{g}^{\left( i \right)}}={{\widetilde{g}}^{\left( i \right)}}+{{C}_{1}}, \quad {{p}^{\left( i \right)}}={{\widetilde{p}}^{\left( i \right)}}-{{C}_{1}} \eqno(30)
$$
is true for~some $i\in \{ 0 \}\cup \mathbb{N}$. According to the~formula~(15), we have the~following representation:
$$
{{g}^{\left( i+1 \right)}}( z )={{\widetilde{g}}^{\left( i+1 \right)}}\left( z \right)+{{C}_{1}}, \quad {{\Phi }_{-}}\left( z \right)\in \left[ {{r}_{i-1}},{{r}_{i}} \right]\!,
$$	
where
$$
{{\tilde{g}}^{\left( i+1 \right)}}( z )={{\mu }_{1}}\left( {{\Phi }_{-}}\left( z \right) \right)-{{\tilde{p}}^{\left( i \right)}}\left( \gamma \left( {{\Phi }_{-}}\left( z \right) \right)+a{{\Phi }_{-}}\left( z \right) \right)-~w\left( {{\Phi }_{-}}\left( z \right),\gamma \left( {{\Phi }_{-}}\left( z \right) \right) \right)+{{C}_{1}}, \quad {{\Phi }_{-}}\left( z \right)\in \left[ {{r}_{i}},{{r}_{i+1}} \right]\!.
$$
In turn, the~formulas (18) and (19) imply that
$$
{{p}^{\left( i+1 \right)}}( z )={{\tilde{p}}^{\left( i+1 \right)}}\left( z \right)+{{C}_{1}}, \quad {{\Phi }_{-}}\left( z \right)\in \left[ l+a{{l}_{i}},l+a{{l}_{i+1}} \right]\!.
$$	
where
\begin{multline*}
    {{\widetilde{p}}^{\left( i+1 \right)}}\!\left( z \right)={{\widetilde{p}}^{\left( i \right)}}\!\left( l+a{{l}_{i}}+0 \right)+ \int\limits_{l+a{{l}_{i}}}^z \exp\! \left( \frac{b\left( l+a{{l}_{i}}-\eta  \right)}{{{a}^{2}}} \right) 
    \bigg( D{{{\widetilde{p}}}^{\left( i+1 \right)}}\!\left( l+a{{l}_{i}}+0 \right) + \\
    +\int\limits_{l+a{{l}_{i}}}^\eta \exp\! \left( \frac{b\left( \xi -l-a{{l}_{i}} \right)}{{{a}^{2}}} \right){{a}^{-2}}\left\{ {{\mu }_{2}}\!\left( \frac{\xi -l}{a} \right)-{{a}^{2}}{{D}^{2}}{{{\widetilde{g}}}^{\left( i \right)}}\!\left( 2l-\xi  \right)-bD{{{\widetilde{g}}}^{\left( i \right)}}\!\left( 2l-\xi  \right)- \right.  \\ 
    \left. -\left. b\frac{\partial w}{\partial x}\!\left( \frac{\xi -l}{a},l \right)-\frac{{{\partial }^{2}}w}{\partial {{t}^{2}}}\!\left( \frac{\xi -l}{a},l \right) \right\}d\xi  \right)d\eta , \quad z\in \left[ l+a{{l}_{i}},l+a{{l}_{i+1}} \right]\!.
\end{multline*}
Thus, the~induction step is proven. Consequently, the~expression (30) holds for~all $i\in \{ 0 \}\cup \mathbb{N}$. Therefore, we can define the~functions $g$ and $p$ as follows 
$$
g=\widetilde{g}+{{C}_{1}}, \quad p=\widetilde{p}-{{C}_{1}},
$$	
where
\begin{align*}
    & \widetilde{g}\!\left( z \right)={{\widetilde{g}}^{\left( 0 \right)}}( z ),\quad z\in [ 0,l ], \\
    & \widetilde{g}\!\left( z \right)={{\widetilde{g}}^{\left( i \right)}}( z ),\quad {{\Phi }_{-}}\!\left( z \right)\in \left[ {{r}_{i-1}},{{r}_{i}} \right]\!,\quad i=1,2,\ldots . \\
    & \widetilde{p}\!\left( z \right)={{\widetilde{p}}^{\left( 0 \right)}}( z ),\quad z\in \left[ 0,l \right]\!, \\
    & \widetilde{p}\!\left( z \right)={{\widetilde{p}}^{\left( j \right)}}( z ),\quad z\in \left[ l+a{{l}_{j-1}},l+a{{l}_{j}} \right]\!, \quad j=1,2,\ldots.
\end{align*}
The~assertion has been proven.

\textbf{Theorem 1.} {\it Let the~smoothness conditions $(23)$ be satisfied. The~mixed problem $(1)$--$(4)$, $(8)$ has a unique solution in the~class ${{C}^{2}}(\overline{Q})$ if and only if the~matching conditions $(25)$--$(28)$ are satisfied. This solution is determined by~the~formulas $(11)$--$(15)$, $(18)$, $(19)$. }

\textbf{Proof.} Assertions 6 and 7 imply the~existence of~a solution. The~uniqueness of~the~solution stems from its construction and Assertion 7 because it is derived from the~general solution. Assuming there are two solutions to the~problem (1)--(4), (8) implies that the~difference between them satisfies the~homogeneous equation (1) and the~homogeneous conditions (3) and (4). Formulas (11)--(15), (18), and (19), as well as Assertion 7, show that the~homogeneous problem has only a zero solution. These results prove the~uniqueness of~the~solution to the~problem (1)--(4).

\textbf{4. Main problem.} Since in the~general case $\psi \notin {{C}^{1}}([0,l])$, the~problem (1)--(4) has no solution belonging to the~class ${{C}^{2}}(\overline{Q})$, i.e., the~problem (1)--(4) has no global classical solution defined on the~set $\overline{Q}$. However, it is possible to define a classical solution on a smaller set $\overline{Q}\setminus \Gamma $ that will satisfy Eq.~(1) on the~set $\overline{Q}\setminus \Gamma $ in the~standard sense and some additional conjugation conditions on the~set $\Gamma $. 

\textbf{Definition 1.} A function $u$ is a \textit{classical solution} of~the~problem (1)--(4) if it is representable in~the~form $u={{u}_{1}}+{{u}_{2}}$, where ${{u}_{1}}$ is a classical solution of~the~problem (1)--(4) with $v=0$ and ${{u}_{2}}$ satisfies Eq. (1) with $f\equiv 0$, the~initial conditions ${{u}_{2}}(0,x)={{\partial }_{t}}{{u}_{2}}(0,x)=0$, $x\in [0,l]$, the~boundary conditions (3) with ${{\mu }_{1}}={{\mu }_{2}}\equiv 0$, and the~following matching conditions
\begin{align}
    & [{{({{u}_{2}})}^{+}}-{{({{u}_{2}})}^{-}}](t,x={{\gamma }_{-}}({{r}_{i}})+at)=0, \quad i=0,1,\ldots, \tag{31} \\
    & [{{({{u}_{2}})}^{+}}-{{({{u}_{2}})}^{-}}](t,x=l+a{{l}_{i}}-at)=0, \quad i=0,1,\ldots, \tag{32} \\
    & \begin{aligned}
        & [{{({{\partial }_{t}}{{u}_{2}})}^{+}}-{{({{\partial }_{t}}{{u}_{2}})}^{-}}](t,x=l+a{{l}_{i}}-at)={{C}^{(i)}},\quad i\in \mathrm{Even}\!\left[ \mathbb{N}\cup \{0\} \right]\!, \\
        & [{{({{\partial }_{t}}{{u}_{2}})}^{+}}-{{({{\partial }_{t}}{{u}_{2}})}^{-}}](t,x=l+a{{l}_{i}}-at)=0,\quad i\in \mathrm{Odd}\!\left[ \mathbb{N} \right]\!,
    \end{aligned} \tag {33}
\end{align}
where ${{C}^{(i)}}$, $i=0,1,...$, are some constants, 
$$\mathrm{Even}\!\left[ \Omega  \right]=\left\{ x : x\in \Omega \wedge x\equiv 0 \pmod 2 \right\}\!,$$
and 
$$\mathrm{Odd}\!\left[ \Omega  \right]=\left\{ x : x\in \Omega \wedge x\equiv 1 \pmod 2 \right\}\!.$$

\textbf{Theorem 2.} {\it Let the~smoothness conditions $(23)$ be satisfied. The~mixed problem $(1)$--$(4)$ has a unique solution in the~sense of~Definition 1 if and only if the~matching conditions $(25)$--$(28)$ are satisfied.} 

\textbf{Proof.} According to Theorem 1, under the~smoothness conditions (23), the~``smooth'' part of~the~solution, i.e., the~function ${{u}_{1}}$ from Definition 1, exists and is unique if and only if the~conditions (25)--(28) are satisfied. The~``discontinuous'' part of~the~solution, i.e., the~function ${{u}_{2}}$ from Definition 1, can be defined by~the~formula
\[{{u}_{2}}(t,x)=g_{*}^{(i)}(x-at)+p_{*}^{(j)}(x+at), \quad (t,x)\in \overline{{{Q}^{(i,j)}}}. \eqno(34)\]
where 	
\begin{align}
    & g_{*}^{(0)}(x)=p_{*}^{(0)}(x)=0, \quad x\in [0,l], \tag{35} \\
    & g_{*}^{(i)}(z)=-p_{*}^{(i-1)}(\gamma ({{\Phi }_{-}}(z))+a{{\Phi }_{-}}(z)), \quad {{\Phi }_{-}}(z)\in [{{r}_{i-1}},{{r}_{i}}], \quad i\in \mathbb{N}, \tag{36} \\
    & \begin{multlined}
        p_{*}^{(j)}(z)=p_{*}^{(j-1)}\left( l+a{{l}_{j-1}}-0 \right) + \\
        +\int\limits_{l+a{{l}_{j-1}}}^z \exp\! \left( \frac{b\left( l+a{{l}_{j-1}}-\eta  \right)}{{{a}^{2}}} \right)\left\{ Dp_{*}^{(j-1)}\left( l+a{{l}_{j-1}}-0 \right)+\frac{{{C}^{(j-1)}}{{\chi }_{\text{Odd}}}(j)}{a}- \right. \\
         - \int\limits_{l+a{{l}_{j-1}}}^\eta {{a}^{-2}}\exp \left( \frac{b\left( \xi -l-a{{l}_{i-1}} \right)}{{{a}^{2}}} \right)\left( {{a}^{2}}{{D}^{2}}g_{*}^{(j-1)}\left( 2l-\xi  \right)+bDg_{*}^{(j-1)}\left( 2l-\xi  \right) \right)d\xi  \bigg\}d\eta,  \\ z\in [l+a{{l}_{j-1}},l+a{{l}_{j}}], \quad j\in \mathbb{N},
    \end{multlined}
    \tag{37}
\end{align}
where ${{\chi }_{\mathrm{Odd}}}$ is an indicator function of~a set $\mathrm{Odd}\!\left[ \mathbb{N} \right]$. The~formulas (34)--(37) can be derived in the~same way as (11)--(15), (18), and (19). Let us write this out in more detail. First, we construct the~function ${{u}_{2}}$ from the~general solution of~Eq. (1) with $f\equiv 0$. It has the~form	
\begin{equation}
    {u}_{2}\!\left( t,x \right)=\tilde{g}_{*}^{\left( i \right)}\!\left( x-at \right)+\tilde{p}_{*}^{\left( j \right)}\!\left( x+at \right), \:\: \left( t,x \right)\in \overline{{{Q}^{\left( i,j \right)}}},\:\: i=0,1,\ldots ,\:\: j=0,1,\ldots, \:\: \left| i-j \right|\le 1.
\tag {38}
\end{equation}
Here, the~functions $\tilde{g}_{*}^{\left( i \right)}$ and $\tilde{p}_{*}^{\left( i \right)}$ ($i=0,1,\ldots )$ are piecewise twice continuously differentiable functions. From the~homogeneous Cauchy conditions (2) with $\varphi \equiv \psi \equiv 0$ it follows [26, p.~174--175]
\begin{gather}
    \tilde{g}_{*}^{\left( 0 \right)}\!\left( x \right)={{\tilde{C}}_{1}},\quad x\in \left[ 0,l \right]\!, \tag{39} \\
    \tilde{p}_{*}^{\left( 0 \right)}\!\left( x \right)=-{{\tilde{C}}_{1}},\quad x\in \left[ 0,l \right]\!, \tag{40}
\end{gather}
where ${{\tilde{C}}_{1}}$ is a real integration constant. The~functions $\tilde{g}_{*}^{\left( i \right)}$ for~all $i\in \mathbb{N}$ are determined from the~Dirichlet boundary condition~(3) with ${{\mu }_{1}}\equiv 0$. For a fixed positive integer $i$ and $j=i-1$, we substitute the~general solution (38) of~Eq. (1) with $f\equiv 0$ into the~boundary condition (3) with ${{\mu }_{1}}\equiv 0$ and obtain
$$
\tilde{g}_{*}^{\left( i \right)}\!\left( \gamma \left( t \right)-at \right)+\tilde{p}_{*}^{\left( i-1 \right)}\!\left( \gamma \left( t \right)+at \right)=0, \quad t\in \left[ {{r}_{i-1}},{{r}_{i}} \right]\!, \quad i=0,1,\ldots. 
$$
Using the~results of~Section 3, we arrive at the~following equality:
$$
\tilde{g}_{*}^{\left( i \right)}\!\left( z \right)=-\tilde{p}_{*}^{\left( i-1 \right)}\!\left( \gamma \left( {{\Phi }_{-}}\!\left( z \right) \right)+a{{\Phi }_{-}}\!\left( z \right) \right)\!, \quad {{\Phi }_{-}}\!\left( z \right)\in \left[ {{r}_{i-1}},{{r}_{i}} \right]\!, \quad i=1,2,\ldots \eqno(41)
$$
Similarly, we define the~function $\tilde{p}_{*}^{\left( j \right)}$ for~all $j\in \mathbb{N}$:
\begin{multline}
    \tilde{p}_{*}^{\left( j \right)}\!\left( z \right)=\tilde{p}_{*}^{\left( j \right)}\!\left( l+a{{l}_{j-1}}+0 \right)+ \int\limits_{l+a{{l}_{j-1}}}^{z}\exp \!\left( \frac{b\left( l+a{{l}_{j-1}}-\eta  \right)}{{{a}^{2}}} \right)\times \\
    \bigg( D\tilde{p}_{*}^{\left( j \right)}\!\left( l+a{{l}_{j-1}}+0 \right)+\underset{l+a{{l}_{j-1}}}{\overset{\eta }{\mathop \int }}\,\exp \left( \frac{b\left( \xi -l-a{{l}_{i-1}} \right)}{{{a}^{2}}} \right){{a}^{-2}}\times  \\
    \times \left. \left\{ {{\mu }_{2}}\!\left( \frac{\xi -l}{a} \right)-{{a}^{2}}{{D}^{2}}\tilde{g}_{*}^{\left( j-1 \right)}\!\left( 2l-\xi  \right)-b D\tilde{g}_{*}^{\left( j-1 \right)}\!\left( 2l-\xi  \right) \right\}d\xi  \right)d\eta , \\ z\in \left[ l+a{{l}_{j-1}},l+a{{l}_{j}} \right]\!,\quad j=1,2,\ldots , \tag{42}
\end{multline}
Here the~values $\tilde{p}_{*}^{\left( j \right)}\!\left( l+a{{l}_{j-1}}+0 \right)$ and $D\tilde{p}_{*}^{\left( j \right)}\!\left( l+a{{l}_{j-1}}+0 \right)$ in the~representation (42) must be chosen so that the~matching conditions (32) and (33) are satisfied. From representation (38) it follows
\begin{gather*}
    \left[ {{\left( {{u}_{2}} \right)}^{+}}-{{\left( {{u}_{2}} \right)}^{-}} \right]\!\left( t,x=l+a{{l}_{j}}-at \right)=\tilde{p}_{*}^{\left( j+1 \right)}\!\left( l+a{{l}_{j}}+0 \right)-\tilde{p}_{*}^{\left( j-1 \right)}\!\left( l+a{{l}_{j}}-0 \right)\!,\quad j=0,1,\ldots, \\
    \begin{multlined}
        \left[ {{\left( \frac{\partial {{u}_{2}}}{\partial t} \right)}^{+}}-{{\left( \frac{\partial {{u}_{2}}}{\partial t} \right)}^{-}} \right]\!\left( t,x=l+a{{l}_{j}}-at \right)=aD\tilde{p}_{*}^{\left( j+1 \right)}\!\left( l+a{{l}_{j}}+0 \right)-aD\tilde{p}_{*}^{\left( j \right)}\!\left( l+a{{l}_{j}}+0 \right)\!,\\ i=0,1,\ldots .
    \end{multlined}
\end{gather*}
So, we need to choose
\begin{align*}
    & \tilde{p}_{*}^{\left( j \right)}\!\left( l+a{{l}_{j-1}}+0 \right)=\tilde{p}_{*}^{\left( j-1 \right)}\!\left( l+a{{l}_{j-1}}-0 \right)\!,~~~~j=1,2,\ldots , \\
    & D\tilde{p}_{*}^{\left( j \right)}\!\left( l+a{{l}_{j-1}}+0 \right)=D\tilde{p}_{*}^{\left( j-1 \right)}\!\left( l+a{{l}_{j-1}}-0 \right)+{{a}^{-1}}{{C}^{(j-1)}}{{\chi }_{\text{Odd}}}\left( j \right)\!,~~~~j=1,2,\ldots .
\end{align*}
Similarly to Assertion 7, it can be shown that the~functions $\tilde{g}_{*}^{\left( j \right)}$ and $p_{*}^{\left( j \right)}$ have the~form
$$
\tilde{g}_{*}^{\left( j \right)}=g_{*}^{\left( j \right)}+{{\tilde{C}}_{1}}, \quad \tilde{p}_{*}^{\left( j \right)}=p_{*}^{\left( j \right)}-{{\tilde{C}}_{1}}, \quad j=0,1,\ldots,
$$
where $g_{*}^{\left( j \right)}$ and $p_{*}^{\left( j \right)}$ are some functions that do not depend on the~constant ${{\tilde{C}}_{1}}$ and are determined by~the~relations (35)--(37). Therefore, since the~expression (38) does not depend on the~constant ${{\tilde{C}}_{1}}$ then the~constructed solution is unique and it is expressed by~the~formulas (38)--(42) when ${{\tilde{C}}_{1}}=0$, i.e., by~the~formulas (35)--(37).

The~fulfillment of~the~conjugation conditions (31)--(33) is verified directly. Let us demonstrate this with the~condition (31) for~$i=1$, i.e., 
$$
\left[ {{\left( {{u}_{2}} \right)}^{+}}-{{\left( {{u}_{2}} \right)}^{-}} \right]\!\left( t,x=at \right)=0. \eqno(43)
$$
The~formulas (34)--(36) imply the~representations  ${{\left( {{u}_{2}} \right)}^{+}}\!\left( t,x=at \right)=0$ and ${{\left( {{u}_{2}} \right)}^{+}}\!\left( t,x=at \right)=0$. These equalities yield $\left[ {{\left( {{u}_{2}} \right)}^{+}}-{{\left( {{u}_{2}} \right)}^{-}} \right]\left( t,x=at \right)=0$. It proves the~condition (43). The~remaining conjugation conditions are verified similarly. The~fact that the~functions $g_{*}^{\left( i \right)}$ and $p_{*}^{\left( i \right)}$, $i=0,1,\ldots $, belong to classes ${{C}^{2}}\!\left( \mathfrak{D}\!\left( g_{*}^{\left( i \right)} \right) \right)$ and ${{C}^{2}}\!\left( \mathfrak{D}\!\left( p_{*}^{\left( i \right)} \right) \right)$, respectively, as established analogously to Assertion 6, allows us to~conclude that the~function ${{u}_{2}}$ belongs to the~class ${{C}^{2}}\!\left( \overline{{{Q}^{\left( i,j \right)}}} \right)$ for~all $i\in \{ 0 \}\cup \mathbb{N}$, $j\in \{ 0 \}\cup \mathbb{N}$, $\left| i-j \right|\le 1,$ and satisfies Eq. (1) on the~sets $\overline{{{Q}^{\left( i,j \right)}}}$ with $f\equiv 0$. The~fulfillment of~the~initial and boundary conditions also follows from the~construction, since the~functions $g_{*}^{\left( i \right)}$ and $p_{*}^{\left( i \right)}$, $i=0,1,\ldots $, are chosen so that the~homogeneous initial and boundary conditions are satisfied. The~theorem is proved.

\textbf{Remark 1.} {\it The~solution to problems $(1)$--$(4)$ is not uniquely defined. Specifically, we must specify the~constants ${{C}^{\left( i \right)}}$ for~all $i\in \mathrm{Even}\!\left[ \mathbb{N} \right]$.}

\textbf{Remark 2.} {\it According to Theorem 2, any choice of~the~constants ${{C}^{\left( i \right)}}$, $i\in \mathrm{Even}\!\left[ \mathbb{N} \right]$, uniquely determines the~solution.}

\textbf{5. Physically correct solution.} By a physically correct solution we mean one that correctly describes the~impact process. The~following statement holds.

\textbf{Assertion 8.} {\it Suppose  $D\gamma \!\left( {{r}_{j}} \right)=0$ for~all $j\in \mathrm{Even}\!\left[ \mathbb{N} \right]$. Then, we can set ${{C}^{\left( j \right)}}=v$ and obtain a physically correct solution in Theorem 2. }

The~proof~of~Assertion 8 is given in the~paper [27].

The~following theorem describes a more general case of~Assertion 8.

\textbf{Theorem 3.} {\it If we set
$$
{{C}^{\left( i \right)}}=v \prod\limits_{j=1}^{i/2} \frac{a+{\gamma }'\!\left( {{r}_{2j-1}} \right)}{a-{\gamma }'\!\left( {{r}_{2j-1}} \right)}, \quad i\in \mathrm{Even}\!\left[ \mathbb{N} \cup \{ 0 \}\right]\!, \eqno(44)
$$
then a solution of~the~problem (1)--(4) constructed in Theorem 2 is physically correct.
}

\textbf{Proof.} The~conditions (30) and (31) are derived from the~continuity. Therefore, we only need to show the~correctness of~the~condition (32) under the~conditions. At the~initial moment $t=0$, the~rod is subjected to an impact at the~end $x=l$. It generates the~shock wave that spreads along the~characteristic $x+at=l$. Its velocity must satisfy the~following condition
$$[{{({{\partial }_{t}}u)}^{+}}-{{({{\partial }_{t}}u)}^{-}}](t,x=l-at)=v.$$
It proves (44) for~$i=0$. For the~derivation of~the~previous equality, we refer the~reader to our paper![28]. 

Furthermore, when the~rod reaches its endpoint, it is immediately reflected and propagates along the~characteristic at a speed that we cannot set but can calculate as follows:
$$
[{{({{\partial }_{t}}u)}^{+}}-{{({{\partial }_{t}}u)}^{-}}](t,x={{\gamma }_{-}}({{r}_{1}})-at)=v\frac{a+{\gamma }'({{r}_{1}})}{a-{\gamma }'({{r}_{1}})}.
$$
Interacting with the~moving end of~the~rod changes its velocity. Since waves propagate at the~same speed in elastic rods~[29], we should have
$$
[{{({{\partial }_{t}}u)}^{+}}-{{({{\partial }_{t}}u)}^{-}}](t,x=l+a{{l}_{2}}-at)=[{{({{\partial }_{t}}u)}^{+}}-{{({{\partial }_{t}}u)}^{-}}](t,x={{\gamma }_{-}}({{r}_{1}})-at)=v\frac{a+{\gamma }'({{r}_{1}})}{a-{\gamma }'({{r}_{1}})}.
$$
It proves (44) for~$i=2$. 

We will prove it further using the~method of~mathematical induction. The~base case has been proven. Now, let us prove the~induction step. Assume that we have
$$
[{{({{\partial }_{t}}u)}^{+}}-{{({{\partial }_{t}}u)}^{-}}](t,x=l+a{{l}_{i}}-at)=v \prod\limits_{j=1}^{i/2} \frac{a+{\gamma }'\!\left( {{r}_{2j-1}} \right)}{a-{\gamma }'\!\left( {{r}_{2j-1}} \right)}
$$
for~some $i\in \mathrm{Even}\!\left[ \mathbb{N} \right]$. The~wave that moves along the~characteristic $l+a{{l}_{i}}-at$ will be reflected from the~end $x=\gamma (t)$ of~the~rod and its speed will become equal to
$$
[{{({{\partial }_{t}}u)}^{+}}-{{({{\partial }_{t}}u)}^{-}}](t,x={{\gamma }_{-}}({{r}_{i+1}})-at)=v\frac{a+{\gamma }'\!\left( {{r}_{i+1}} \right)}{a-{\gamma }'\!\left( {{r}_{i+1}} \right)} \prod\limits_{j=1}^{i/2} \frac{a+{\gamma }'\!\left( {{r}_{2j-1}} \right)}{a-{\gamma }'\!\left( {{r}_{2j-1}} \right)}=v\prod\limits_{j=1}^{1+ i/2}\frac{a+{\gamma }'\!\left( {{r}_{2j-1}} \right)}{a-{\gamma }'\!\left( {{r}_{2j-1}} \right)}
$$
according to the~formulas (34) – (37). The~wave propagates at the~same speed in elastic rods, so 
$$
[{{({{\partial }_{t}}u)}^{+}}-{{({{\partial }_{t}}u)}^{-}}](t,x=l+a{{l}_{i+2}}-at)=[{{({{\partial }_{t}}u)}^{+}}-{{({{\partial }_{t}}u)}^{-}}](t,x={{\gamma }_{-}}({{r}_{i+1}})-at)=v\prod\limits_{j=1}^{1+i/2}\frac{a+{\gamma }'\!\left( {{r}_{2j-1}} \right)}{a-{\gamma }'\!\left( {{r}_{2j-1}} \right)}.
$$
The~induction step is proven. Thus, we have proven the~physical correctness of~the~condition (33) for~all $i\in \mathrm{Even}\!\left[ \mathbb{N}\cup \{0\} \right]$ when the~equalities~(44) are satisfied. Since there were initially no shock waves in the~rod, no shock wave propagates along the~characteristics 
$$
x=at, \quad x=l+a{{l}_{1}}-at, \quad x={{\gamma }_{-}}({{r}_{2}})+at, \quad x=l+a{{l}_{2}}-at, \quad x={{\gamma }_{-}}({{r}_{4}})+at 
$$
etc., i. e., $x={{\gamma }_{-}}({{r}_{2i}})+at$ and $x=l+a{{l}_{2i+1}}-at$, $i=0,1,...$. This proves the~physical correctness of~the~condition (33) for~all $i\in \mathrm{Odd}\!\left[ \mathbb{N} \right]$. 

\textbf{Conclusions.} In the~present paper, we have obtained the~necessary and sufficient conditions under which a unique classical solution of~a mixed problem exists for~the~wave equation with discontinuous conditions in a curvilinear half-strip. We have constructed the~solution in an implicit analytical form. We have proposed a method for~constructing solutions to mixed problems for~hyperbolic equations with discontinuous conditions in curvilinear domains.

\begin{center}{\textbf{References}}\end{center}

\small{1. Kil'chevskii N. A. \textit{Teoriya soudarenii tverdykh tel} (Theory of~Collisions between Solid Bodies), Kiev, Naukova Dumka, 1969.}

\small{2. Boussinesq J. Du choc longitudinal d'une barre prismatique, fix\'{e}e \`{a} un bout et heurt\'{e}e \`{a} l'autre, \textit{Comptes Rendus}, 1883, vol.~97, pp.~154--157.}

\small{3. Koshlyakov N. S., Smirnov M. M., Gliner E. B. \textit{Differential Equations of~Mathematical Physics}, Amsterdam, North-Holland Publishing Co., 1964.}

\small{4. Nikolai E. L. On the~theory of~longitudinal impact of~elastic rods, \textit{Trudy Leningr. industr. in-ta}, 1939, no. 3, p. 85--93.}

\small{5. Manzhosov V. K. \textit{Modeli prodolnogo udara} (Longitudinal Impact Models), Ulyanovsk, Ulyanovsk State Technical University, 2006.}

\small{6. Bityurin A. A., Manzhosov V. K. \textit{Prodolnyi udar neodnorodnogo sterzhnya o zhestkuyu pregradu} (Longitudinal Impact of~a Non-Uniform Rod on a Rigid Barrier), Ulyanovsk, Ulyanovsk State Technical University, 2009.}

\small{7. Slepukhin V. V. Modeling of~Wave Processes under Longitudinal Impact in Rod Systems of~Homogeneous Structure, \textit{Abstract of~Cand. Sci. (Techn.) Dissertation}, Ulyanovsk, 2009.}

\small{8. Zhilin P. A. \textit{Prikladnaya mekhanika: teoriya tonkikh uprugikh sterzhnei} (Applied Mechanics: Theory of~Thin Elastic Rods), Saint Petersburg, Polytechnic University Publishing House, 2007.}

\small{9. Gaiduk S. I. Some problems related to the~theory of~longitudinal impact on a rod, \textit{Differential Equations}, 1977, vol.~12(1976), pp.~607--617.}

\small{10. Rasulov M. L. \textit{Methods of~Contour Integration}, Amsterdam, North-Holland Publishing Co., 1967.}

\small{11. Korzyuk V. I., Rudzko J. V. A mathematical investigation of~one problem of~the~longitudinal impact on an elastic rod with an elastic attachment at the~end. \textit{Tr. Inst. Mat., Minsk}, 2023, vol.~31, no. 1, pp.~81--87.}

\small{12. Yufeng X., Dechao Z. Analytical solutions of~impact problems of~rod structures with springs, \textit{Comput. Methods Appl. Mech. Eng.}, 1998, vol.~160, pp.~315--323.}

\small{13. Hu B., Eberhard P. Symbolic computation of~longitudinal impact waves, \textit{Computer Methods in Applied Mechanics and Engineering}, 2001, vol.~190, no. 37-38, pp.~4805--4815.}

\small{14. Etiwa R. M., Elabsy H. M., Elkaranshawy H. A. Dynamics of~longitudinal impact in uniform and composite rods with effects of~various support conditions, \textit{Alexandria Engineering Journal}, 2023, vol.~65, pp.~1–22.}

\small{15. Gomez B. J., Repetto C. E., Stia C. R., Welti R. Oscillations of~a string with concentrated masses, \textit{Eur. J. Phys.}, 2007, vol.~28, pp.~961--975.}

\small{16. Tikhonov A. N., Samarskii A. A. \textit{Equations of~Mathematical Physics}, New York, Dover Publications, 2011.}

\small{17. Naumavets S. N. Classical solution of~the~first mixed problem for~the~one-dimensional wave equation with a differential polynomial of~the~second order in the~boundary conditions. \textit{Zbirnyk statei. Matematyka. Informatsiini tekhnologi\"{\i}. Osvita}, 2018, no. 5, P. 96--101.}

\small{18. Kapustin N. Yu. On spectral problems arising in the~theory of~the~parabolic-hyperbolic heat equation. \textit{Dokl. Math.}, 1996, vol.~54, pp.~607--610.}

\small{19. Kapustin N. Yu., Moiseev E. I. Spectral problems with the~spectral parameter in the~boundary condition. \textit{Differ. Equ.}, 1997, vol.~33, no. 1, pp.~116--120.}

\small{20. Korzyuk V. I., Kozlovskaya I. S., Naumavets S. N. Classical Solution of~the~First Mixed Problem for~the~Wave Equation in a Curvilinear Half-Strip, \textit{Differential Equations}, 2020, vol.~56, pp.~98--108.}

\small{21. Korzyuk V. I., Stolyarchuk I. I. Classical solution of~the~first mixed problem for~second-order hyperbolic equation in curvilinear half-strip with variable coefficients, \textit{Differential Equations}, 2017, vol.~53, pp.~74--85.}

\small{22. Korzyuk V. I., Rudzko J. V. Classical Solution of~the~First Mixed Problem for~the~Telegraph Equation with a Nonlinear Potential in a Curvilinear Quadrant, \textit{Differential Equations}, 2023, vol.~59, pp.~1075--1089.}

\small{23. Ammari K., Bchatnia A., El Mufti K. Stabilization of~the~wave equation with moving boundary, \textit{European Journal of~Control}, 2018, vol.~39, pp.~35--38.}

\small{24. Liu L., Gao H. The~stabilization of~wave equations with moving boundary, \textit{arXiv:2103.13631}.}

\small{25. De Jesus I. P., Lapa E. C., Limaco J. Controllability for~the~wave equation with moving boundary, \textit{Electronic Journal of~Differential Equations}, 2021, vol.~2021, no. 60, pp.~1–12.}

\small{26. Korzyuk V. I. \textit{Uravneniya matematicheskoi fiziki} (Equations of~Mathematical Physics), Moscow, URSS, 2021.}

\small{27. Korzyuk V. I., Rudzko J. V., Kolyachko V. V. Classical solution of~a mixed problem for~the~wave equation with discontinuous initial conditions in a curvilinear half-strip, \textit{Doklady of~the~National Academy of~Sciences of~Belarus}, 2025 (in print).}

\small{28. Korzyuk V. I., Rudzko J. V., Kolyachko V. V. Solutions of~problems with discontinuous conditions for~the~wave equation, \textit{Journal of~the~Belarusian State University. Mathematics and Informatics}, 2023, vol.~3, pp.~6--18.}

\small{29. Zhuravkov M., Lyu Y., Starovoitov E. \textit{Mechanics of~Solid Deformable Body}, Singapore, Springer, 2023.}

\end{document}